\documentclass[a4paper,11pt,oneside, reqno]{amsart}

\usepackage{amsaddr}
\usepackage{cite}
\usepackage{hyperref}
\usepackage{lmodern} 
\usepackage{bbm}
\usepackage[bb=boondox]{mathalfa}

\usepackage{tikz}

\newcommand{\one}{\mathbb{1}}
\newcommand{\ZZ}{\mathbb{Z}}
\newcommand{\RR}{\mathbb{R}}
\newcommand{\CC}{\mathbb{C}}
\newcommand{\NN}{\mathbb{N}}
\newcommand{\cO}{\mathcal{O}}
\newcommand{\SSS}{\mathbb{S}}
\newcommand{\eps}{\varepsilon}
\newcommand{\rmi}{\mathrm{i}}
\newcommand{\mn}{\mathrm{min}}

\DeclareMathOperator{\supp}{\mathrm{supp}}

\DeclareMathOperator{\dom}{\mathrm{dom}}
\DeclareMathOperator{\spec}{\mathrm{spec}}

\DeclareMathOperator{\spece}{\mathrm{spec}_\mathrm{ess}}
\DeclareMathOperator{\sgn}{\mathrm{sgn}}

\newtheorem{theorem}{Theorem}
\newtheorem{corollary}[theorem]{Corollary}

\renewcommand{\Tilde}{\widetilde}

\newcommand{\dd}{\mathrm{d}}

\begin{document}
	\title[]{\large A self-adjoint indefinite Laplacian with~a~large~essential spectrum}

	\author[]{Konstantin Pankrashkin}
	\address{\footnotesize Carl von Ossietzky Universit\"at Oldenburg,
		Fakult\"at V -- Mathematik und Naturwissenschaften,  Institut f\"ur Mathematik, Ammerl\"ander Heerstr. 114--118, 26129~Oldenburg (Oldb), Germany}
	\email{konstantin.pankrashkin@uol.de}

	\date{}

\begin{abstract}
We study a self-adjoint operator $L=-\mathrm{div}\,h\,\mathrm{grad}$ in a bounded domain with a piecewise constant function $h$ changing
its sign along a hypersurface. Such operators arise as
 mathematical models of negative metamaterials, and earlier works found that
the essential spectrum of $L$ can be non-empty. However, in all known examples the non-trivial essential spectrum
consisted of a single point. We show that the geometry can be chosen in such a way that the
essential spectrum of $L$ covers the whole real axis. This is achieved by a direct
construction of quasimodes that are localized near the set on which the sign change of $h$ occurs.
\end{abstract}

\maketitle

\section{Introduction}

\subsection{Motivation}
The present note may be viewed as a continuation of the works~\cite{bk,cpp} that highlighted
the appearance of a non-trivial essential spectrum for indefinite Laplacians in bounded domains.

Let $\Omega\subset\RR^n$ be a bounded domain and $\Sigma$ a hypersurface cutting $\Omega$ into two parts $\Omega_\pm$. For $\mu>0$, consider the function $h$ whose value is $1$ in $\Omega_-$ and $(-\mu)$ in $\Omega_+$. We are interested in self-adjoint differential operators $L$ in $L^2(\Omega)$ formally given by the differential expression
\[
L:\quad u\mapsto -\nabla\cdot (h\nabla u)
\]
(usually equipped with Dirichlet boundary condition on $\partial\Omega$). Such operators arise in the mathematical theory of negative-index metamaterials: they model an interaction of the region $\Omega_-$ filled by a dielectric with the region $\Omega_+$ filled by a metamaterial, and the two regions touch and interact along the \emph{interface} $\Sigma$. The parameter $\mu$ is often referred to as a \emph{contrast}. For the physical background we refer
to~\cite{ves}, and for a survey of mathematical results on negative-index materials, including applications to superlensing and cloaking, to the monograph~\cite{amm} and the review~\cite{ngu17}.  For a more rigorous definition one uses the natural identification
\[
L^2(\Omega)\simeq L^2(\Omega_-)\oplus L^2(\Omega_+), \quad 
u\simeq(u_-,u_+),
\quad
u_\pm:= u|_{\Omega_\pm},
\]
then the above operator $L$ acts as $(u_-,u_+)\mapsto (-\Delta u_-,\mu\Delta u_+)$ with $u_\pm$
satisfying the boundary and transmission conditions
\begin{equation}
	\label{transm0}
u_\pm=0 \text{ on } \partial\Omega\cap\partial\Omega_\pm, 
\quad u_-=u_+ \text{ on } \Sigma,
\quad \partial_- u_-= \mu\partial_+ u_+ \text{ on } \Sigma,
\end{equation}
where $\partial_\pm$ are the outer (with respect to $\Omega_\pm$) normal derivatives on $\Sigma$, and we denote the functions and their various Sobolev traces on the boundary by the same symbols for the sake of readability. Finding the precise regularity of $u_\pm$ guaranteeing the self-adjointness of $L$ turns out to be an involved problem depending on the geometric properties and smoothness of~$\Sigma$. If $\Omega_\pm$ are Lipschitz, it is natural to start with the ``minimal'' realization $L_\mn$,
\begin{equation}
	\label{eq-ll}
\begin{aligned}
L_\mn:&\ (u_-,u_+)\mapsto (-\Delta u_-,\mu\Delta u_+),\\
\dom L_\mn&:=\Big\{
(u_-,u_+)\in H^2(\Omega_-)\times H^2(\Omega_+): \ \eqref{transm0} \text{ is satisfied}
\Big\},
\end{aligned}
\end{equation}
which is easily shown to be a symmetric operator in $L^2(\Omega)$. If $\mu\ne 1$ (non-critical contrast), both $\Sigma$ and $\partial\Omega$ are smooth, and $\Sigma$ is closed and  does not touch the exterior boundary $\partial\Omega$,
then it is known that $L_\mn$ is self-adjoint \cite{bbdr,cpp}. The same conclusion holds for $\mu\ne 1$ and some choices
of~$\Sigma$ that touch the exterior boundary, see e.g. \cite[Sec.~5]{hussein} or~\cite[Sec.~8]{ss2}. 
Some self-adjointness results for non-smooth $\Sigma$ with a range of $\mu\ne 1$ are available as well, see \cite{bbdr,bccj,bp-rmc}, and curious transition phenomena between smooth and non-smooth cases are known~\cite{chen}. In all these cases (with $\mu\ne 1$), the operator $L_\mn$ turns out to have compact resolvent, so its spectrum consists of isolated eigenvalues of finite multiplicities that accumulate at $\pm\infty$, while the positive and negative eigenvalues satisfy a Weyl-like bound \cite[Sec.~6.2]{hussein}. The more recent work \cite{cdv} obtains similar conclusions in a more detailed form for an analog of the above operator on manifolds. The papers \cite{cz,cdv} noted the existence of quasimodes localized near the interface.

It was already noted in \cite{bbdr} that the case $\mu=1$ (critical contrast) is expected to be significantly different. While some earlier announcements appeared in \cite{kostr}, it seems that the self-adjointness problem was rigorously addressed in \cite{bk} by looking at the two-dimensional
case with $\Omega_\pm$ being squares sharing a common side $\Sigma$. In this case,
the above operator $L_\mn$ is not self-adjoint, but has a natural self-adjoint extension which
turns out to have an infinitely degenerate zero eigenvalue, i.e. the essential spectrum becomes non-empty.
In \cite{cpp}, one looked more systematically at the case when $\Sigma$ is a smooth closed submanifold that encloses $\Omega_-$ and does not touch the exterior boundary $\partial\Omega$. Then it was shown \cite[Thm.~1]{cpp} that the above operator $L_\mn$ (with $\mu=1$) is not closed but essentially self-adjoint, while the domain and the spectral properties
of its unique self-adjoint extension $L$ depend on the dimension and the geometry of $\Sigma$. In particular, by \cite[Thm.~2]{cpp}
for $n=2$ one always has $\spece L=\{0\}$, and $\dom L$ is not contained in any $H^s(\Omega\setminus\Sigma)$ for $s>0$.
On the other hand, if $n\ge 3$ and the principal curvatures of $\Sigma$
are either all strictly positive or all strictly negative  on each connected component of $\Sigma$,
then $\dom L\subset H^1(\Omega\setminus\Sigma)$,
and $L$ has compact resolvent, see \cite[Thm.~3]{cpp}. If $n\ge 3$ and a portion  of $\Sigma$ is flat, then $\{0\}\subset\spece L$.
 These observations
were based on the reduction of the analysis of $L$ to the study of a pseudodifferential operator on $\Sigma$, which
is smoothing for $n=2$ and has order one with curvature-dependent principal
symbol for $n\ge 3$ (in particular, elliptic, if the above sign condition for the curvatures is satisfied). Later papers \cite{faikl,kroum} provided further examples (featuring a kind of separation of variables) of self-adjoint operators $L$ of the above type with $\mu=1$ and $\spece L=\{0\}$. We note that the essential spectrum of indefinite Laplacians
is closely related to cloaking effects in metamaterials, see \cite{bouch,kett,kohn,ngu16} or \cite[Ch.~10-11]{amm} for various approaches.

Given the above observations, it is natural to ask whether the essential spectrum of a self-adjoint indefinite Laplacian with critical contrast ($\mu=1$) can be larger than just a point. As an additional motivation one may mention the parallel activity on Dirac operators with  transmission conditions, for which it was recently found that critical couplings in bounded domains may produce curvature-dependent intervals of essential spectrum in three dimensions, while in two dimensions the non-trivial essential spectrum is a finite set~\cite{BBHP}. These results can be seen as a localization of the preceding
spectral studies of Dirac operators for the whole space \cite{BHOP20,BP24}.

\subsection{Main result}

Our goal is to show that the spectral effect of the critical contrast $\mu=1$
can be significantly more striking than everything that was observed in the previous works: 
it is possible to choose the bounded domain $\Omega$ and the hypersurface $\Sigma$ in the above constructions in such a way that the essential spectrum
of the respective  indefinite Laplacian $L$ covers the whole real axis. In other words, no non-trivial spectral enclosure is possible for general geometries. It is particularly surprising that this peculiar spectral effect is realized by a very simple geometry: a single cylindrical portion of the interface $\Sigma$ is sufficient to produce a half-line of essential spectrum, and two cylindrical portions with opposite orientations yield the whole real axis.

{}From now on we consider the three-dimensional case $n=3$ and restrict the geometric setting. Let $\Omega\subset\RR^3$ be a bounded smooth domain and $\Omega_-\subset\Omega$ be a bounded smooth domain
with $\overline{\Omega_-}\subset\Omega$. Denote
\[
\Sigma:=\partial\Omega_-,\quad \Omega_+:=\Omega\setminus\overline{\Omega_-}.
\]
Let $\mu=1$. Then the operator $L_\mn$ in \eqref{eq-ll} has a unique self-adjoint extension $L$ due to \cite[Thm.~1]{cpp}, which 
we call critical indefinite Laplacian.

\begin{theorem}\label{thm1}
Assume that for some $(x_0,y_0)\in \RR^2$ and $a,b,c,d,R\in \RR$ with $0<a<R<b$ and $c<d$ one has
the inclusion
\[
\omega:=\big\{(x,y,z)\in\RR^3:\  a^2<(x-x_0)^2+(y-y_0)^2<b^2,\ c<z<d\big\}\subset \Omega
\]
and the equalities
\begin{align}
\Sigma\cap \omega&=\big\{(x,y,z)\in\omega:\  (x-x_0)^2+(y-y_0)^2=R^2\big\},\nonumber\\
\Omega_-\cap \omega&=\big\{(x,y,z)\in\omega:\  (x-x_0)^2+(y-y_0)^2<R^2\big\},\label{omplus1}\\
\Omega_+\cap \omega&=\big\{(x,y,z)\in\omega:\  (x-x_0)^2+(y-y_0)^2>R^2\big\}.\label{omplus2}
\end{align}
Then for any $E\in(-\infty,0]$ one can explicitly construct a sequence $(u_m)_{m\in\NN}\subset \dom L\setminus\{0\}$ such that $\supp u_m\subset\omega$ and
\begin{equation}
	\dfrac{\|(L-E)u_m\|_{L^2(\Omega)}}{\|u_m\|_{L^2(\Omega)}}\xrightarrow{m\to\infty}0. \label{leu}
\end{equation}
In particular, it holds that $(-\infty,0]\subset\spece L$.
\end{theorem}

We note that the existence of the Weyl sequence $(u_m)$ with \eqref{leu} implies $E\in\spec L$. As this holds for any $E\le0$
and each non-isolated point of the spectrum belongs to the essential spectrum, one obtains the last conclusion.

By applying the sign flip one can interchange the roles of $\Omega_+$ and $\Omega_-$:
\begin{corollary}\label{cor1}
Assume that in Theorem~\ref{thm1} the conditions \eqref{omplus1}--\eqref{omplus2}
are replaced with
\begin{align*}
\Omega_+\cap \omega&=\big\{(x,y,z)\in\omega:\  (x-x_0)^2+(y-y_0)^2<R^2\big\},\\
\Omega_-\cap \omega&=\big\{(x,y,z)\in\omega:\  (x-x_0)^2+(y-y_0)^2>R^2\big\}.
\end{align*}
Then
$[0,\infty)\subset\spece L$.	
\end{corollary}

By combining both previous results we obtain:

\begin{corollary}\label{cor2}
	Assume that for some $(x_k,y_k)\in \RR^2$ and $a_k,b_k,c_k,d_k,R_k\in \RR$ with $0<a_k<R_k<b_k$ and $c_k<d_k$, $k\in\{1,2\}$, one has
	the inclusions
	\[
	\omega_k:=\big\{(x,y,z)\in\RR^3:\  a^2_k<(x-x_k)^2+(y-y_k)^2<b_k^2,\ c_k<z<d_k\big\}\subset \Omega,\quad k\in\{1,2\},
	\]
	and the equalities
	\begin{align*}
		\Sigma\cap \omega_k&=\big\{(x,y,z)\in\omega_k: (x-x_k)^2+(y-y_k)^2=R_k^2\big\},\ k\in\{1,2\},\\[\smallskipamount]
		\Omega_-\cap \omega_1&=\big\{(x,y,z)\in\omega_1: (x-x_1)^2+(y-y_1)^2<R^2_1\big\},\\
		\Omega_+\cap \omega_1&=\big\{(x,y,z)\in\omega_1: (x-x_1)^2+(y-y_1)^2>R^2_1\big\},\\[\smallskipamount]
		\Omega_-\cap \omega_2&=\big\{(x,y,z)\in\omega_2: (x-x_2)^2+(y-y_2)^2>R^2_2\big\},\\
\Omega_+\cap \omega_2&=\big\{(x,y,z)\in\omega_2: (x-x_2)^2+(y-y_2)^2<R^2_2\big\}.	
	\end{align*}
	Then $\spec L=\spece L=\RR$.	
\end{corollary}	

\begin{proof}
The conditions for $k=1$ imply $(-\infty,0]\subset \spece L$ by Theorem~\ref{thm1}, and the conditions for $k=2$ yield
$[0,\infty)\subset\spece L$ by Corollary~\ref{cor1}. It remains to recall that $\spece L\subset\spec L\subset\RR$ due to the self-adjointness.
\end{proof}

To be more specific, pick $0<r<R$ and $\ell>0$, and let $\Omega_-$ be any bounded smooth domain
that coincides with the cylindrical layer $\{r^2<x^2+y^2<R^2\}$
in the slab $\{|z|<\ell\}$, and let $\Omega$ be any bounded smooth domain containing $\overline{\Omega_-}$.
Then the essential spectrum of the respective self-adjoint critical indefinite Laplacian $L$ coincides with $\RR$ by Corollary~\ref{cor2}.

The rest of the text is devoted to the proof of Theorem~\ref{thm1}. The construction of the Weyl sequence  $(u_m)$ is fully elementary: the functions $u_m$  are given by \eqref{um-vm} and \eqref{um1} with the coordinate change \eqref{eq-Phi}, and they are exponentially localized near $\Sigma$ for $m\to\infty$. The main effort is spent on the estimate of the respective norms with full details.

\section{Proof of Theorem~\protect\ref{thm1}}

\subsection{Change of variables}\label{sec21}
Without loss of generality we may assume $(x_0,y_0)=(0,0)$, $R=1$ and $c=-d$ with $d>0$ (otherwise one applies suitable dilations and translations, which does not influence the sought spectral inclusion).
For any function $u\in\dom L$ supported in $\omega$ we have
\[
(L u)(x,y,z)=\sgn(x^2+y^2-1)\Delta u(x,y,z) \text{ in }\omega\setminus\Sigma.
\]

Recall that for a function $u$ defined in $\omega$ and $w(r,\theta,z):=u(r\cos\theta,r\sin\theta,z)$ one has
\[
\iiint |u|^2 \dd x\,\dd y\,\dd z=\iiint |w|^2r\,\dd r\,\dd\theta\,\dd z
\]
(with corresponding integration domains) and
\[
(\Delta u)(r\cos\theta,r\sin\theta,z)=\big(\partial_r^2 w+r^{-1}\partial_r w+r^{-2}\partial_\theta^2 w+\partial_z^2 w\big)(r,\theta,z).
\]
Further, for $v(t,\theta,z):=w(e^t,\theta,z)$ we have
\[
\iiint |u|^2 \dd x\,\dd y\,\dd z=\iiint |v|^2e^{2t}\,\dd t\,\dd\theta\,\dd z
\]
(with corresponding integration domains) as well as
\begin{gather*}
(\partial_r^2 w+r^{-1}\partial_r w)(e^t,\theta,z)=e^{-2t}\partial^2_t v(t,\theta,z),\\
(\Delta u)(e^t\cos\theta,e^t\sin\theta,z)=\Big(e^{-2t}(\partial^2_t v+\partial^2_\theta v)+\partial^2_z v\Big)(t,\theta,z).
\end{gather*}

Being motivated by the above observations, pick a small $\eps\in (0,d)$ such that
$e^{2\eps}\le 2$ and the map
\begin{equation}
	\label{eq-Phi}
\Phi:\   (t,\theta,z)\mapsto (e^t\cos\theta,e^t\sin\theta,z)
\end{equation}
defines a diffeomorphism between
\[
\Theta:=(-\eps,\eps)\times\SSS\times (-\eps,\eps),\qquad \SSS:=\RR/(2\pi\ZZ),
\]
and an open subset $\Phi(\Theta)$ of $\omega$. For any function $u\in L^2(\Omega)$ supported in $\Phi(\Theta)$, by denoting
\[ 
v:=u\circ \Phi
\]
we have
\begin{gather*}
\tfrac{1}{2}\|v\|^2_{L^2(\Theta)}\le	\|u\|^2_{L^2(\Omega)}\equiv \iiint_\Theta |v|^2 e^{2t}\,\dd t\,\dd\theta\,\dd z\le 2\|v\|^2_{L^2(\Theta)},\quad
\|v\|^2_{L^2(\Theta)}:=\iiint_\Theta |v|^2 \dd t\,\dd\theta\,\dd z.
\end{gather*}	
If additionally $u\in \dom L$, then
\begin{equation*}
	(L u)\big(\Phi(t,\theta,z)\big)=(\sgn t)\Big(e^{-2t}(\partial^2_t v+\partial^2_\theta v)+\partial^2_z v\Big)(t,\theta,z) \text{ in } \Theta\cap \{t\ne 0\}.
\end{equation*}	
We also remark that the above functions $u$ are smooth on each side of $\Sigma$ if and only if the respective functions $v$ are smooth
on each side of the surface $\{t=0\}\subset\Theta$, and in this case the transmission conditions \eqref{transm0} for $u$ are satisfied (with $\mu=1$) if and only if
\begin{equation}
	\label{transm1}
	v(0^+,\theta,z)=v(0^-,\theta,z),\quad \partial_t v(0^-,\theta,z)+\partial_t v(0^+,\theta,z)=0,\quad (\theta,z)\in \SSS\times(-\eps,\eps).
\end{equation}

In view of the above computations, to arrive at a sequence $(u_m)$ as in \eqref{leu}
with required localization properties it is sufficient to construct
a sequence of functions $v_m$ that are compactly supported in $\Theta$, smooth up to the boundary in $\{\pm t>0\}\cap\Theta$, satisfy the transmission conditions \eqref{transm1} for each $m$ and such that
\begin{gather*}
	\dfrac{\|\Tilde v_m-E v_m\|_{L^2(\Theta)}}{\|v_m\|_{L^2(\Theta)}}\xrightarrow{m\to\infty}0,\\
	\Tilde v_m(t,\theta,z):=(\sgn t)\Big[e^{-2t}(\partial^2_t v_m+\partial^2_\theta v_m)+\partial^2_z v_m\Big](t,\theta,z),\quad t\ne 0,
\end{gather*}	
then one can simply set
\begin{equation}
	\label{um-vm}
u_m:=v_m\circ \Phi^{-1} \text{ extended by zero to the whole of $\Omega$.}
\end{equation}
Indeed, by construction the functions $u_m$ are supported in $\Phi(\Theta)\subset\omega$, smooth on both sides of $\Sigma$
and satisfy the transmission conditions~\eqref{transm0} on $\Sigma$, so $u_m\in\dom L_\mn\subset\dom L$.
We also have $\Tilde v_m-Ev_m= \big[(L-E)u_m\big] \circ\Phi$, therefore,
\begin{gather*}
\|(L-E)u_m\|^2_{L^2(\Omega)}\le 2\|\Tilde v_m-E v_m\|^2_{L^2(\Theta)},\quad
\|u_m\|^2_{L^2(\Omega)}\ge\tfrac{1}{2}\|v_m\|^2_{L^2(\Theta)},\\
\bigg(\dfrac{\|(L-E)u_m\|_{L^2(\Omega)}}{\|u_m\|_{L^2(\Omega)}}\bigg)^2\le 4 \bigg(\dfrac{\|\Tilde v_m-E v_m\|_{L^2(\Theta)}}{\|v_m\|_{L^2(\Theta)}}\bigg)^2\xrightarrow{m\to\infty}0.
\end{gather*}
In the rest of the text we construct the functions $v_m$ with the above properties.

\subsection{Candidates for {\boldmath $v_m$}}

Let $E\le 0$, which will be fixed for the rest of the text, and
pick $\chi\in C^\infty_c(-\eps,\eps)$ with $\|\chi\|_{L^2(-\eps,\eps)}=1$.
For $m\in\NN$ consider the function
\[
f_m: \ z\mapsto \chi(z)e^{\rmi\sqrt{|mE|}z},
\]
then $\|f_m\|_{L^2(-\eps,\eps)}=1$. For any $k\in\NN_0$ and $m\in\NN$ we have
\[
\Big|f_m^{(2k)}(z)\Big|=\bigg|\sum_{j=0}^{2k}{\textstyle\binom{2k}{j}} \chi^{(2k-j)}(z) \rmi^j|mE|^\frac{j}{2}e^{\rmi\sqrt{|mE|}z}\bigg|\le
m^k\sum_{j=0}^{2k}2^{2k}|E|^\frac{j}{2}\Big| \chi^{(2k-j)}(z)\Big|,
\]
which implies
\begin{equation}
	\label{fnorm}
\|f_m^{(2k)}\|_{L^2(-\eps,\eps)}=\cO(m^k) \text{ for }m\to\infty.
\end{equation}
For the function $r_m$ defined by
\begin{align}
r_m(z):=\dfrac{1}{m}f^{(2)}_m(z)-E f_m(z)&=\dfrac{2\rmi\chi'(z)\sqrt{|mE|}e^{\rmi \sqrt{|Em|}z} +\chi''(z)e^{\rmi \sqrt{|mE|}z}}{m}\label{bm1}\\
&\equiv 2\rmi \sqrt{|E|}\dfrac{\chi'(z)e^{\rmi \sqrt{|mE|}z}}{\sqrt{m}}+\dfrac{\chi''(z)e^{\rmi \sqrt{|mE|}z}}{m}\nonumber
\end{align}
we have
\begin{equation}
	\label{bm2}
\|r_m\|^2_{L^2(-\eps,\eps)}=\cO(m^{-1}),\quad m\to\infty.
\end{equation}

Let us pick any $\delta\in(0,\eps)$ and $\eta\in C^\infty_c(-\eps,\eps)$ with $\one_{(-\delta,\delta)}\le\eta\le \one_{(-\eps,\eps)}$. For $m\in\NN$, define
functions $v_m$ on $\Theta$ by
\begin{align}
	v_m(t,\theta,z)&:=\eta(t)e^{\rmi m\theta}e^{-m|t|}g_m(t,z), \label{um1}\\
	g_m(t,z)&:=
	 f_m(z)+\dfrac{|t|+t|t|}{2m}f^{(2)}_m(z)
	 + \dfrac{mt^2+|t|}{8m^3}f^{(4)}_m(z).\nonumber
\end{align}
We show that $v_m$ the properties stated in Subsection~\ref{sec21}.

\subsection{Transmission conditions for {\boldmath$v_m$}}

The functions $v_m$ are smooth in $\{\pm t\ge0\}\cap\Theta$ and continuous on $\Theta$,
\[
	v_m(0^+,\theta,z)=v_m(0^-,\theta,z),\quad (\theta,z)\in\SSS\times(-\eps,\eps),
\]
i.e. the first condition in \eqref{transm1} is satisfied.
Further, 
\begin{align*}
	\partial_t v_m(t,\theta,z)&=\eta'(t)e^{\rmi m\theta}e^{-m|t|}g_m(t,z)
		+\eta(t)e^{\rmi m\theta}e^{-m|t|}\Big( -m(\sgn t)g_m(t,z) + \partial_t g_m(t,z)\Big).
\end{align*}
Noting that $\eta'=0$ in $(-\delta,\delta)$ we obtain
\[
	\partial_t v_m(0^+,\theta,z)+\partial_t v_m(0^-,\theta,z)= \eta(0)e^{\rmi m\theta}\Big(\partial_t g_m(0^+,z)+
	\partial_t g_m(0^-,z)\Big).
\]
Due to
\begin{align*}
	\partial_t g_m(t,z)&=
\dfrac{\sgn t +2|t|}{2m}f^{(2)}_m(z)
+ \dfrac{2mt+\sgn t}{8m^3}f^{(4)}_m(z),\\
	\partial_t g_m(0^\pm,z)&=\pm\Big( \dfrac{1}{2m}f^{(2)}_m(z)+\dfrac{1}{8m^3}f^{(4)}_m(z)\Big)
\end{align*}
it follows that $\partial_t v_m(0^+,\theta,z)+\partial_t v_m(0^-,\theta,z)=0$,
which is the second condition in \eqref{transm1}.

\subsection{Lower bound for the norm of {\boldmath$v_m$}}
For any $x,y,z\in\CC$ we have
\begin{align*}
	|x+y+z|^2&\ge |x|^2-2|x||y+z|+|y+z|^2\ge |x|^2-\Big(\dfrac{1}{2}|x|^2 + 2|y+z|^2\Big)+|y+z|^2\\
	&=\dfrac{1}{2}|x|^2-|y+z|^2\ge \dfrac{1}{2}|x|^2-2\big(|y|^2+|z|^2\big).
\end{align*}
By applying this elementary bound to the expression for $v_m$ we obtain, with a suitable $c>0$ independent of $m$,
\begin{align*}
	\big|v_m(t,\theta,z)\big|^2&\ge\dfrac{1}{2}\eta(t)^2e^{-2m|t|}\big|f_m(z)\big|^2\\
	&\quad-2\eta(t)^2e^{-2m|t|}\bigg( \dfrac{(t+1)^2t^2}{4m^2}\big| f^{(2)}_m(z)\big|^2 + \dfrac{(m|t|+1)^2t^2}{64m^6}\big| f^{(4)}_m(z)\big|^2\bigg)\\
	&\ge\dfrac{1}{2}\one_{(-\delta,\delta)}(t)e^{-2m|t|}\big|f_m(z)\big|^2-c\one_{(-\eps,\eps)}(t)e^{-2m|t|}\bigg( \dfrac{t^2}{m^2}\big| f^{(2)}_m(z)\big|^2 + \dfrac{t^2}{m^4}\big| f^{(4)}_m(z)\big|^2\bigg),
\end{align*}
hence,
\begin{multline*}
	\|v_m\|^2_{L^2(\Theta)}=\int_{-\eps}^{\eps}\int_0^{2\pi}\int_{-\eps}^{\eps} \big|v_m(t,\theta,z)\big|^2\dd t\,\dd\theta\,\dd z
	\ge \pi \int_{-\delta}^{\delta}e^{-2m|t|} \dd t\big\|f_m\big\|^2_{L^2(-\eps,\eps)}\\
	{}-2\pi c\bigg[\dfrac{1}{m^2}\int_{-\eps}^{\eps}|t|^2 e^{-2m|t|}\dd t \big\|f_m^{(2)}\big\|^2_{L^2(-\eps,\eps)}
	+\dfrac{1}{m^4}\int_{-\eps}^{\eps}t^2e^{-2m|t|}\dd t \big\|f_m^{(4)}\big\|^2_{L^2(-\eps,\eps)}\bigg].
\end{multline*}	
Recall that $\|f_m\|^2_{L^2(-\eps,\eps)}=1$. For any $k\in\NN_0$ we have
\begin{equation}
	\label{intmk}
\begin{aligned}
	\int_{-\eps}^\eps |t|^k e^{-2m|t|}\dd t&=2\int_0^\eps t^k e^{-2mt}\dd t\le 2\int_0^\infty t^k e^{-2mt}\dd t\\
	&=\dfrac{2}{(2m)^{k+1}}\int_0^\infty y^k e^{-y}\dd y=\dfrac{\Gamma(k+1)}{2^k m^{k+1}}\equiv \dfrac{k!}{2^k m^{k+1}},
\end{aligned}
\end{equation}
while $\|f_m^{(2k)}\|^2_{L^2(-\eps,\eps)}\le c_k m^{2k}$ with $c_k>0$ by \eqref{fnorm}. It follows that
\begin{equation}
	\label{vmlow}
\begin{aligned}
	\|v_m\|^2_{L^2(\Theta)}&\ge \dfrac{\pi}{m}(1-e^{-2m\delta})-2\pi c\Big(
	\dfrac{1}{m^2}\cdot\dfrac{2}{4m^3}\cdot c_1 m^2 + \dfrac{1}{m^4}\cdot\dfrac{2}{4m^3} \cdot c_2 m^4
	\Big)\\
	&= \dfrac{\pi}{m}(1-e^{-2m\delta})-2\pi c
	\dfrac{c_1+c_2}{2m^3}\ge \dfrac{\pi}{2 m}
\end{aligned}
\end{equation}
for all sufficiently large $m$.

\subsection{Computing the error {\boldmath $\widetilde{v_m}-E v_m$}}

In all computations below $t\ne 0$. Recall that
\begin{align*}
\Tilde v_m(t,\theta,z)=(\sgn t)\Big[e^{-2t}(\partial^2_t v_m +\partial^2_{\theta} v_m) +\partial_z^2 v_m\Big](t,\theta,z).
\end{align*}

Note that
\begin{align}
	\partial^2_t g_m(t,z)&=
	\dfrac{\sgn t }{m}f^{(2)}_m(z)
	+ \dfrac{1}{4m^2}f^{(4)}_m(z), \nonumber\\
	\partial_z^2g_m(t,z)&=
	f_m^{(2)}(z)+\dfrac{|t|+t|t|}{2m}f^{(4)}_m(z)
	+ \dfrac{mt^2+|t|}{8m^3}f^{(6)}_m(z), \label{dz2g}
\end{align}
and compute
\begin{align}
\partial^2_{\theta} v_m(t,\theta,z)&=\eta(t)e^{\rmi m\theta}e^{-m|t|}\Big[-m^2 g_m(t,z)\Big],\nonumber\\
\partial_t^2 v_m(t,\theta,z)&=
\eta(t)e^{\rmi m\theta}e^{-m|t|}\Big[ m^2 g_m(t,z)-2m(\sgn t) \partial_t g_m(t,z)+\partial^2_t g_m(t,z) \Big]+w_m(t,\theta,z),\nonumber\\
w_m(t,\theta,z)&:=\eta''(t)e^{\rmi m\theta}e^{-m|t|}g_m(t,z)\nonumber\\
&\qquad+2\eta'(t)e^{\rmi m\theta}e^{-m|t|}\Big[-m(\sgn t) g_m(t,z)+\partial_tg_m(t,z)\Big], \label{eq-wm}
\end{align}
hence,
\begin{align*}
	(\partial^2_t v_m +\partial^2_{\theta} v_m)(t,\theta,z)&=\eta(t)e^{\rmi m\theta}e^{-m|t|}\Tilde g_m(t,z)+w_m(t,\theta,z),\\
	\Tilde g_m(t,z)&:=-2m(\sgn t) \partial_t g_m(t,z)+\partial^2_t g_m(t,z)\\
	&=-2m(\sgn t)\Big( \dfrac{\sgn t+2|t|}{2m}f^{(2)}_m(z)
	+ \dfrac{2mt+\sgn t}{8m^3}f^{(4)}_m(z)\Big)\\
	&\quad +\Big( \dfrac{\sgn t }{m}f^{(2)}_m(z)
	+ \dfrac{1}{4m^2}f^{(4)}_m(z)\Big)\\
	&=-\Big( (1+2t)f^{(2)}_m(z)
	+ \dfrac{2m|t|+1}{4m^2}f^{(4)}_m(z)\Big) +\Big( \dfrac{\sgn t}{m}f^{(2)}_m(z)
	+ \dfrac{1}{4m^2}f^{(4)}_m(z)\Big)\\
	&=\Big( -1-2t + \dfrac{\sgn t}{m}\Big)f^{(2)}_m(z)-\dfrac{|t|}{2m}f^{(4)}_m(z).
\end{align*}	

By representing $e^{-2t}=1-2t+t^2b(t)$ with $b\in L^\infty(-\eps,\eps)$
we obtain, with suitable functions $b_j\in L^\infty(-\eps,\eps)$ indepedent of $m$,
\begin{gather*}
	e^{-2t}(\partial^2_t v_m +\partial^2_{\theta} v_m)(t,\theta,z)=\eta(t)e^{\rmi m\theta}e^{-m|t|} h_m(t,z)+b_1(t)w_m(t,\theta,z),\\
	h_m(t,z):=\Big( -1 + \dfrac{\sgn t}{m} +t^2b_2(t)+\dfrac{t b_3(t)}{m}\Big)f^{(2)}_m(z)-\dfrac{|t|+t^2b_5(t)}{2m}f^{(4)}_m(z),
\end{gather*}	
and with \eqref{dz2g} it follows that
\begin{align}
	\Tilde v_m(t,\theta,z)&=
	(\sgn t)\Big[\eta(t)e^{\rmi m\theta}e^{-m|t|}p_m(t,z)+b_1(t)w_m(t,\theta,z)\Big], \label{del1}\\
	p_m(t,z)&:=h_m(t,z)+\partial_z^2g_m(t,z)\nonumber\\
	\label{pm1}
	&=\dfrac{\sgn t}{m}f^{(2)}_m(z) +\Big(t^2b_2(t) +\dfrac{t b_3(t)}{m}\Big)f^{(2)}_m(z)
+\dfrac{t^2 b_6(t)}{m}f^{(4)}_m(z)+\dfrac{mt^2+|t|}{8m^3}f^{(6)}_m(z).
\end{align}
Using~\eqref{bm1} to rewrite the first summand in \eqref{pm1} to obtain
\begin{align*}
	p_m(t,z)&=(\sgn t)\Big(E f_m(z)+r_m(z)\Big)+\Big(t^2b_2(t) +\dfrac{t b_3(t)}{m}\Big)f^{(2)}_m(z)\\
	&\qquad+\dfrac{t^2 b_6(t)}{m}f^{(4)}_m(z)+\dfrac{mt^2+|t|}{8m^3}f^{(6)}_m(z)\\
	&=(\sgn t)Ef_m(z)+q_m(t,z),\\
	q_m(t,z)&:=(\sgn t)r_m(z)+\Big(t^2b_2(t) +\dfrac{t b_3(t)}{m}\Big)f^{(2)}_m(z)+\dfrac{t^2 b_6(t)}{m}f^{(4)}_m(z)+\dfrac{mt^2+|t|}{8m^3}f^{(6)}_m(z).
\end{align*}

The substitution into \eqref{del1} gives
\begin{align*}
\Tilde v_m(t,\theta,z)&=	E \eta(t)e^{\rmi m\theta}e^{-m|t|}f_m(z)+(\sgn t)\eta(t)e^{\rmi m\theta}e^{-m|t|}q_m(t,z)+(\sgn t)b_1(t)w_m(t,\theta,z),
\end{align*}
and by \eqref{um1} the first summand on the right-hand side can be represented as
\[
E \eta(t)e^{\rmi m\theta}e^{-m|t|}f_m(z)=E v_m(t,\theta,z)
-E\eta(t)e^{\rmi m\theta}e^{-m|t|}\Big(\dfrac{|t|+|t|t}{2m}f^{(2)}_m(z)
+\dfrac{mt^2+|t|}{8m^3}f^{(4)}_m(z)\Big),
\]
which results in
\begin{gather}
	(\Tilde v_m -E v_m)(t,\theta,z)=\rho_m(t,\theta,z)+(\sgn t)b_1(t)w_m(t,\theta,z),\label{err1}\\
	\begin{aligned}
		\rho_m(t,\theta,z):=\eta(t)e^{\rmi m\theta}e^{-m|t|}\Big[&{}-E\dfrac{|t|+|t|t}{2m}f^{(2)}_m(z)
	-E \dfrac{mt^2+|t|}{8m^3}f^{(4)}_m(z)+(\sgn t)q_m(t,z)\Big].
	\end{aligned}\label{rhom}
\end{gather}

\subsection{The norm of the error}

With suitable $C_j>0$ independent of $m$ we estimate
\begin{align*}
	\big|q_m(t,z)\big|^2&\le C_1\bigg[ \big|r_m(z)\big|^2 +\Big(t^4+\dfrac{t^2}{m^2}\Big)\big|f^{(2)}_m(z)\big|^2+\dfrac{t^4}{m^2}|f^{(4)}_m(z)\big|^2+\Big(\dfrac{t^4}{m^4}+\dfrac{t^2}{m^6}\Big)|f^{(6)}_m(z)\big|^2\bigg],\\
	\big|\rho_m(t,\theta,z)\big|^2&\le C_2 e^{-2m|t|}\bigg[ \dfrac{t^2}{m^2}\big|f^{(2)}_m(z)\big|^2+\Big(\dfrac{t^4}{m^4}+\dfrac{t^2}{m^6}\Big)\big|f^{(4)}_m(z)\big|^2+\big|q_m(t,z)\big|^2
	\bigg]\\
	&\le C_3 e^{-2m|t|}\bigg[ \big|r_m(z)\big|^2+ \Big(t^4+\dfrac{t^2}{m^2}\Big)\big|f^{(2)}_m(z)\big|^2+\Big(\dfrac{t^4}{m^2}+\dfrac{t^2}{m^6}\Big)\big|f^{(4)}_m(z)\big|^2\\
	&\qquad\qquad+
	\Big(\dfrac{t^4}{m^4}+\dfrac{t^2}{m^6}\Big)\big|f^{(6)}_m(z)\big|^2
	\bigg].
\end{align*}

By \eqref{rhom} it follows that
\begin{equation}
	\label{rhom1}
\begin{aligned}
	\|\rho_m\|_{L^2(\Theta)}^2&\le C_4\int_{-\eps}^{\eps}e^{-2m|t|}\bigg[ \|r_m\|^2_{L^2(-\eps,\eps)}+ \Big(t^4+\dfrac{t^2}{m^2}\Big)\big\|f^{(2)}_m\big\|^2_{L^2(-\eps,\eps)}\\
	&\qquad +\Big(\dfrac{t^4}{m^2}+\dfrac{t^2}{m^6}\Big)\big\|f^{(4)}_m\big\|^2_{L^2(-\eps,\eps)}+
	\Big(\dfrac{t^4}{m^4}+\dfrac{t^2}{m^6}\Big)\big\|f^{(6)}_m\big\|^2_{L^2(-\eps,\eps)}\bigg]\,\dd t\\
	&
	\stackrel{\text{\eqref{fnorm},\eqref{bm2}}}\le C_5\int_{-\eps}^{\eps}e^{-2m|t|}\bigg[ m^{-1}+ \Big(t^4+\dfrac{t^2}{m^2}\Big)m^2+\Big(\dfrac{t^4}{m^2}+\dfrac{t^2}{m^6}\Big)\cdot m^4\\
	&\qquad+
	\Big(\dfrac{t^4}{m^4}+\dfrac{t^2}{m^6}\Big)m^6\bigg]\,\dd t\\
	&\le C_6\int_{-\eps}^{\eps}e^{-2m|t|}\bigg[
	m^{-1}+t^2+m^2t^4\bigg]\,\dd t\\
	&\stackrel{\eqref{intmk}}{\le} C_7\Big( m^{-1}\cdot m^{-1}+m^{-3}+m^2 m^{-5}\Big)\le C_8 m^{-2}.
\end{aligned}	
\end{equation}

Using the expression \eqref{eq-wm} for $w_m$ and the fact that $\eta'$ and $\eta''$ vanish on $(-\delta,\delta)$ we estimate, with a suitably chosen $\kappa>0$ and for all sufficiently large $m$,
\begin{align}
	\big|w_m(t,\theta,z)\big|^2&\le C_9m^\kappa\one_{\delta<|t|<\eps}e^{-2m|t|}\sum_{k=0}^2\big|f^{(2k)}_m(z)\big|^2
	\le C_9 m^\kappa e^{-2m\delta}\sum_{k=0}^2 \big|f^{(2k)}_m(z)\big|^2, \nonumber\\
\|w_m\|^2_{L^2(\Theta)}&=\int_{-\eps}^\eps \int_\SSS\int_{-\eps}^\eps	\big|w_m(t,\theta,z)\big|^2\dd z\,\dd \theta\,\dd t
\le 4\pi\eps C_9 m^\kappa e^{-2m\delta}\sum_{k=0}^2 \|f^{(2k)}_m\|^2_{L^2(-\eps,\eps)}\nonumber\\
&\le 4\pi\eps C_9 m^\kappa e^{-2m\delta}\sum_{k=0}^2 c_k m^{2k}\le C_{10} e^{-m\delta}.\label{bd-wm}
\end{align}

\subsection{Conclusion}

By \eqref{err1}, for $m\to\infty$ we have the upper error bound
\begin{align*}
	\|\Tilde v_m-E v_m\|^2_{L^2(\Theta)}&\le 2\|\rho_m\|^2_{L^2(\Theta)}+2\|b_1\|^2_{L^\infty(-\eps,\eps)}\|w_m\|^2_{L^2(\Theta)}
	\stackrel{\eqref{rhom1},\eqref{bd-wm}}{=}\cO(m^{-2}),
\end{align*}
as well as the lower bound $\|v_m\|^2_{L^2(\Theta)}\ge \tfrac{\pi}{2m}$ from \eqref{vmlow},
which yields the sought behavior
\[
\dfrac{\|\Tilde v_m -E v_m\|^2_{L^2(\Theta)}}{\|v_m\|^2_{L^2(\Theta)}}=\cO\Big(\dfrac{1}{m}\Big),\quad m\to\infty.
\]

\subsection*{Acknowledgments}
The author was partially supported by the Deutsche Forschungsgemeinschaft (DFG, German Research Foundation) -- 491606144.

\subsection*{Use of artificial intelligence} The quasimode construction was deduced through an interaction with ChatGPT by testing several geometric settings and analytic ans\"atze. The final computation was verified and written up by the author.

\subsection*{Data availability}

No datasets were generated or analyzed during the current study.

\subsection*{Competing interests}

The author has no competing interests to declare that are relevant to the contents of this article.

\end{document}